# PROPERTIES AND REFINEMENTS OF THE FUSED LASSO

By Alessandro Rinaldo[1]

*Carnegie Mellon University*

We consider estimating an unknown signal, both blocky and sparse, which is corrupted by additive noise. We study three interrelated least squares procedures and their asymptotic properties. The first procedure is the fused lasso, put forward by Friedman et al. [*Ann. Appl. Statist.* **1** (2007) 302–332], which we modify into a different estimator, called the fused adaptive lasso, with better properties. The other two estimators we discuss solve least squares problems on sieves; one constrains the maximal $\ell_1$ norm and the maximal total variation seminorm, and the other restricts the number of blocks and the number of nonzero coordinates of the signal. We derive conditions for the recovery of the true block partition and the true sparsity patterns by the fused lasso and the fused adaptive lasso, and we derive convergence rates for the sieve estimators, explicitly in terms of the constraining parameters.

**1. Introduction.** We consider the nonparametric regression model

$$y_i = \mu_i^0 + \varepsilon_i, \qquad i = 1, \ldots, n,$$

where $\mu^0 \in \mathbb{R}^n$ is the unknown vector of mean values to be estimated using the observations $y$, and the errors $\varepsilon_i$ are assumed to be independent with either Gaussian or sub-Gaussian distributions and bounded variances. We are concerned with the more specialized settings where $\mu^0$ can be both *sparse*, with a possibly very large number of zero entries, and *blocky*, meaning that the number of coordinates where $\mu^0$ changes its values can be much smaller than $n$. Figure 1 shows an instance of data generated by corrupting a blocky and sparse signal with additive noise (see Section 2.4 for details about this example). Formally, we assume that there exists a partition $\{\mathcal{B}_1^0, \ldots, \mathcal{B}_{J_0}^0\}$

Received May 2008; revised August 2008.

[1]Supported in part by NSF Grant DMS-06-31589 and a grant from the Pennsylvania Department of Health through the Commonwealth Universal Research Enhancement Program.

**Signal plus noise**

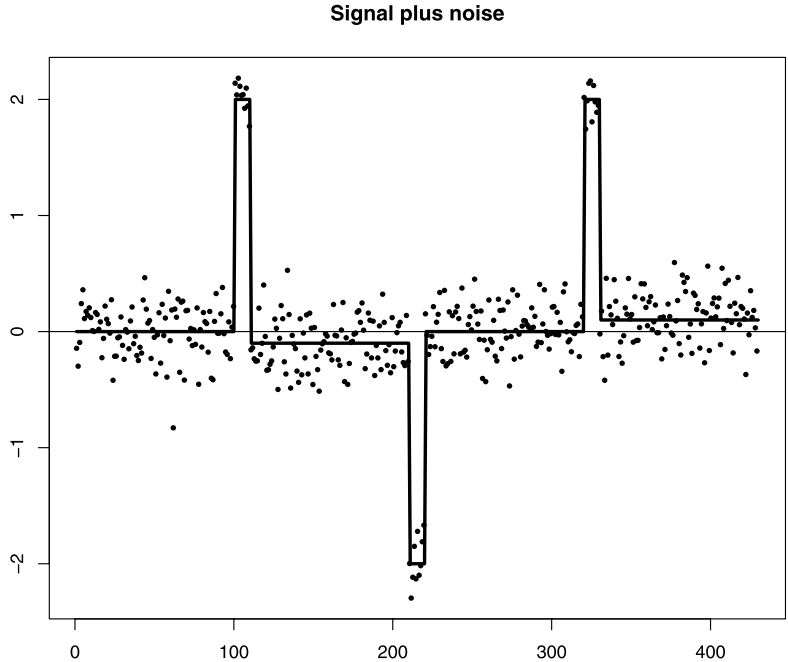

Fig. 1. *Signal (solid line) plus noise for the example described in Section 2.4.*

of $\{1,\ldots,n\}$ into sets of consecutive indexes, from now on called a *block partition*, and a vector $\nu^0 \in \mathbb{R}^{J_0}$, which may be sparse, such that the true mean vector can be written as

$$\mu^0 = \sum_{j=1}^{J_0} \nu_j^0 1_{\mathcal{B}_j^0}, \tag{1.1}$$

where $1_\mathcal{B}$ is the indicator function of the set $\mathcal{B} \subseteq \{1,\ldots,n\}$ (i.e., the $n$-dimensional vector whose $i$th coordinate is 1 if $i \in \mathcal{B}$ and 0 otherwise). The partition $\{\mathcal{B}_1^0, \ldots, \mathcal{B}_{J_0}^0\}$, its size $J_0$, the vector $\nu^0$ of block values and its zero coordinates are all unknown, and our goal is to produce estimates of those or related quantities that are accurate when $n$ is large enough.

In particular, we investigate the asymptotic properties of three different but interrelated methods for the recovery of the unknown mean vector $\mu^0$ under the assumption (1.1).

The first methodology we study, which is the central focus of this work, is the *fused lasso* procedure of Friedman et al. (2007). The *fused lasso* is the penalized least squares estimator

$$\widehat{\mu}^{\mathrm{FL}} = \operatorname*{arg\,min}_{\mu \in \mathbb{R}^n} \left\{ \sum_{i=1}^{n}(y_i - \mu_i)^2 + 2\lambda_{1,n}\|\mu\|_1 + 2\lambda_{2,n}\|\mu\|_{\mathrm{TV}} \right\}, \tag{1.2}$$

**Fusion Adaptive Lasso Estimate**

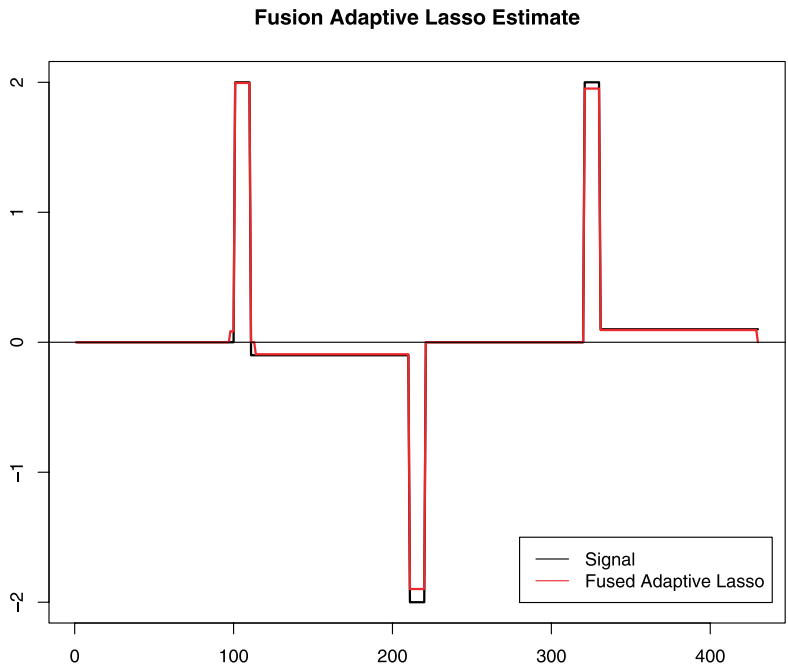

FIG. 2. *A fusion adaptive lasso estimate for the example from Section 2.4, using the most biased fusion estimator shown in Figure 3 the oracle threshold for the lasso penalty, as described in Section 2.3.*

where $\|\mu\|_1 \equiv \sum_{i=1}^n |\mu_i|$ is the $\ell_1^n$ norm and $\|\mu\|_{\text{TV}} \equiv \sum_{i=2}^n |\mu_i - \mu_{i-1}|$ the total variation seminorm of $\mu$, respectively, and $(\lambda_{1,n}, \lambda_{2,n})$ are positive tuning parameters to be chosen appropriately. The solution to the convex program (1.2) can be computed in a fast and efficient way using the algorithm developed in Friedman et al. (2007), where the properties of the fused lasso solution are considered from the optimization theory standpoint. Our analysis has led us to propose a modified version of the fused lasso, which we call the fused adaptive lasso, that has improved properties. Figure 2 shows an example of a fused adaptive lasso fit to the the data displayed in Figure 1.

In our second approach, we turn to a different convex optimization program, namely

$$\operatorname*{arg\,min}_{\mu \in \mathbb{R}^n} \sum_{i=1}^n (y_i - \mu_i)^2$$
(1.3)
$$\text{s.t. } \|\mu\|_1 \leq L_n, \|\mu\|_{\text{TV}} \leq T_n$$

for some nonnegative constants $L_n$ and $T_n$. Notice that, in this alternative formulation, which is akin to the least squares method on sieves, a solution different from $y$ is obtained provided that $\|y\|_1 > L_n$ or $\|y\|_{\mathrm{TV}} > T_n$. The link with the fused lasso estimator is clear. The objective function in the fused lasso problem (1.2) is the Lagrangian function of (1.3), and, in fact, the two problems are equivalent from the point of view of convexity theory.

Our third and final method for the recovery of a sparse and blocky signal is also related to sieve least square procedures, and is more naturally tailored to the model assumption (1.1). Specifically, we study the solution to the highly nonconvex optimization problem

$$\underset{\mu \in \mathbb{R}^n}{\arg\min} \sum_{i=1}^n (y_i - \mu_i)^2$$

(1.4)

$$\text{s.t. } |\{i : \mu_i \neq 0\}| \leq S_n, 1 + |\{i : \mu_i - \mu_{i-1} \neq 0, 2 \leq i \leq n\}| \leq J_n,$$

where $S_n$ and $J_n$ are nonnegative constants. Although lack of convexity makes this problem computationally difficult when $n$ is large, the theoretical relevance of this third formulation stems from the fact that (1.3) is, effectively, a convex relaxation of (1.4).

Our approach to the study of the estimators defined by (1.2), (1.3) and (1.4) is asymptotic, as we allow the block representation for the unobserved signal $\mu^0$ to change with $n$ in such a way that the recovery of a noisy signal under the model (1.1) may become increasingly difficult. Despite being quite closely related as optimization problems, from an inferential perspective, the three procedures under investigation each shed some light on different and, in some way, complementary aspects of this problem.

Overall, our analysis yields conditions for consistency of the block partition and block sparsity estimates by model (1.2) and its variant described in Section 2.3, and explicit rates of consistency of both sieve solutions (1.3) and (1.4). In essence, our results provide conditions for the sequences of regularization parameters $\lambda_{1,n}$, $\lambda_{2,n}$, $L_n$, $J_n$ and $S_n$ to guarantee various degrees of recovery of $\mu^0$.

The article is organized as follows. In Section 2, we study the fused lasso estimator. After deriving an explicit formula for the fused lasso solution in Section 2.1, we establish conditions under which the fused lasso procedure is both sparsistent, in the sense of being a weakly consistent estimator of the partitions, and of the set of nonzero coordinates of $\mu^0$. In Section 2.3, we propose a simple modification of the fused lasso, which we call the fused adaptive lasso, that achieves sparsistency under milder conditions and also allows us to derive an oracle inequality for the empirical risk. Finally, in Section 3, we derive consistency rates for the estimators defined in (1.3) and (1.4), which depend explicitly on the parameters $L_n$ and $T_n$, and of $S_n$ and $J_n$, respectively. The proofs are relegated to the Appendix.

We conclude this introductory section by fixing the notation that we will be using throughout the article. For a vector $\mu \in \mathbb{R}^n$, we let $\mathcal{S}(\mu) = \{i : \mu_i \neq$ 0} denote its support and $\mathcal{J}(\mu) = \{i : \mu_i = \mu_{i-1} \neq 0, i \geq 2\}$ the set of coordinates where $\mu$ changes its value. Furthermore, notice that we can always write

$$\mu = \sum_{j=1}^{J} \nu_j 1_{\mathcal{B}_j}$$

from some (possibly trivial) block partition $\{\mathcal{B}_1, \ldots, \mathcal{B}_J\}$, with $1 \leq J \leq n$, and some vector $\nu \in \mathbb{R}^J$. Then, we will write $\mathcal{JS}(\mu) = \{j : \nu_j \neq 0\}$ for the sets of nonzero blocks of $\mu$. On a final note, although all the quantities defined so far may change with $n$, for ease of readability, we do not always make this dependence explicit in our notation.

1.1. *Previous works and comparison.* The idea of using the total variation seminorm in penalized least squares problem has been exploited and studied in many applications (e.g., signal processing, parametric regression, nonparametric regression and image denoising). From the algorithmic viewpoint, this idea was originally brought up by Rudin, Osher and Fatemi (1992) [for more recent developments, see, e.g., Dobson and Vogel (1997) and Caselles, Chambolle and Novaga (2007), and also DeVore (1998)]. The original motivation for this article was the recent work by Friedman et al. (2007), who devise efficient coordinate-wise descent algorithms for a variety of convex problems. In particular, they propose a novel approach based on a penalized least squares problems using simultaneously the total variation and the $\ell_1$ penalties, which favors solutions that are both blocky and sparse. In the classical nonparametric framework of function estimation, two important contributions in the development and analysis of total variation-based methods come from Mammen and van de Geer (1997) and Davies and Kovac (2001a). Specifically, Mammen and van de Geer (1997) show that least squares splines with adaptively chosen knots are solutions to nonparametric least squares penalized regression problems with total variation penalties and derive, among other things, consistency rates for both the one- and two-dimensional case. Using a different approach, Davies and Kovac (2001a) devise a very simple and effective procedure with $O(n)$ complexity, called the taut-string algorithm, which effectively solves least squares problems with total variation penalty. The taut-string can be used to consistently estimate at an almost optimal rate the number and location of local maxima of an unknown function on $[0, 1]$. Both methods impose very little assumptions on the degree of smoothness of the true underlying function. More recently, Boysen et al. (2009) study jump-penalized least squares regression problems, where the underlying function is assumed to be a linear combination of indicator functions of intervals in $[0, 1]$, and derive consistency rates under different metrics on functional spaces.

Our work differs from the contributions based on a nonparametric function estimation framework of, in particular, Mammen and van de Geer (1997) and Davies and Kovac (2001a) in various aspects, some of which are closely related to the methodology and scope of Friedman et al. (2007). First and foremost, we are interested in the asymptotic recovery of the coordinates of the mean vector $\mu^0$ under the model assumption (1.1), and do not necessarily view them as $n$ evaluations of an unknown function defined on $[0, 1]$. Secondly, we explicitly impose a double asymptotic framework in which the model complexity and the features of the underlying signal change with $n$. This, in particular, allows us to include cases in our analysis where the number of blocks or the number of local extremes grow unbounded with $n$, a feature which typically cannot be directly accommodated in the nonparametric framework. Nonetheless, we remark that there is a simple reformulation of our problem as nonparametric function estimation one. In fact, suppose that we observe $n$ datapoints of the form

$$y_i = n \int_{(i-1)/n}^{i/n} \mu^0(t) \, dt + \varepsilon_i, \qquad i = 1, \ldots, n,$$

from an unknown function $\mu^0 : [0, 1] \to \mathbb{R}$. Setting $\mu_i^0 = n \int_{(i-1)/n}^{i/n} \mu^0(t) \, dt$ would return our original model [see also Boysen et al. (2009) for a similar model]. Furthermore, for the analysis of Section 2, we are only concerned with the simultaneous recovery of both the block partition and of the sparsity pattern of $\mu^0$ and virtually ignore any other features of the signal. On the one hand, this allows us to derive rather strong results, namely sparsistency and the oracle inequality of Theorem 2.7. On the other hand, those results are truly meaningful only when our modeling assumptions (1.1) of a blocky and sparse signal hold, and our analysis should not be expected to be robust to misspecification. In particular, the fused lasso and adaptive fused lasso algorithms should not be expected to work well, both in practice and in theory, with different kinds of signals.

**2. Properties and refinements of the fused lasso estimator.** The crucial feature of the fused lasso solution (1.2), which makes it ideal for the present problem, is that it is simultaneously blocky, because of the total variation penalty $\|\cdot\|_{\mathrm{TV}}$, and sparse, because of the $\ell_1$ penalty $\|\cdot\|_1$. The central goal of this section is to characterize the asymptotic behavior of the regularization parameters $\lambda_{1,n}$ and $\lambda_{2,n}$, so that, as $n \to \infty$, the blockiness and sparsity pattern of the the fused lasso estimates match the ones of the unknown signal $\mu^0$ with overwhelming probability. We first consider the fused lasso estimator as originally proposed in Friedman et al. (2007) and then a simple variant, the fused adaptive lasso, which has better asymptotic properties. For this modified version, we also derive an oracle inequality. We will make the following simplifying assumption on the errors:

(E) The errors $\varepsilon_i$, $1 \leq i \leq n$ are identically distributed centered Gaussian variables with variance $\sigma_n^2$ such that $\sigma_n \to 0$.

In the typical scenario we have in mind, $\sigma_n = \frac{\sigma}{\sqrt{n}}$. Assumption (E) is by no means necessary, and it can be easily relaxed to the case of sub-Gaussian errors.

2.1. *The fused lasso solution.* Below, we provide an explicit formula for the fused lasso solution that offers some insight on its properties and suggests possible improvements. By inspecting (1.2), as both penalty functions $\|\cdot\|_1$ and $\|\cdot\|_{\mathrm{TV}}$ are convex and the objective function is strictly convex, $\widehat{\mu}^{\mathrm{FL}}$ is uniquely determined as the solution to the subgradient equation

$$\widehat{\mu}^{\mathrm{FL}} = y - \lambda_{1,n} s_1 - \lambda_{2,n} s_2, \tag{2.1}$$

where $s_1 \in \partial \|\widehat{\mu}^{\mathrm{FL}}\|_1$ and $s_2 \in \partial \|\widehat{\mu}^{\mathrm{FL}}\|_{\mathrm{TV}}$. For a vector $x \in \mathbb{R}^n$, the subgradient $\partial \|x\|_1$ is a subset of $\mathbb{R}^n$ consisting of vectors $s$ such that $s_i = \mathrm{sgn}(x_i)$, where, with some abuse of notation, we will denote with $\mathrm{sgn}(\cdot)$ the (possibly set-valued) function on $\mathbb{R}$ given by

$$\mathrm{sgn}(x) = \begin{cases} 1, & \text{if } x > 0, \\ -1, & \text{if } x < 0, \\ z, & \text{if } x = 0, \end{cases}$$

where $z$ is any number in $[-1, 1]$. The subgradient $\partial \|x\|_{\mathrm{TV}}$ has a slightly more elaborated form, which is given in Lemma A.1 in the Appendix.

An explicit expression for $\widehat{\mu}^{\mathrm{FL}}$ can be obtained in terms of the *fusion* estimator

$$\widehat{\mu}^F = \underset{\mu \in \mathbb{R}^n}{\arg\min} \left\{ \sum_{i=1}^n (y_i - \mu_i)^2 + 2\lambda_2 \|\mu\|_{\mathrm{TV}} \right\}. \tag{2.2}$$

Notice that, by the same arguments used above, $\widehat{\mu}^F$ is also unique. This fusion estimator solves a regularized least squares problem with a penalty on the total variation of the signal and works by fusing together adjacent coordinates that have similar values to produce a blocky estimate of the form (1.1). We remark that, in the nonparametric function estimation settings, one can obtain $\widehat{\mu}^F$ as a piecewise-constant variable-knot spline function on $[0,1]$ [see Mammen and van de Geer (1997), Proposition 8] and that the taut-string algorithm of Davies and Kovac (2001a) solves the constrained version of (2.2).

For a given solution $\widehat{\mu}^F$ to (2.2), there exists a block partition $\{\widehat{\mathcal{B}}_1, \ldots, \widehat{\mathcal{B}}_{\widehat{J}}\}$ and a unique vector $\widehat{\nu} \in \mathbb{R}^{\widehat{J}}$ such that

$$\widehat{\mu}^F = \sum_{j=1}^{\widehat{J}} \widehat{\nu}_j 1_{\widehat{\mathcal{B}}_j}. \tag{2.3}$$

We take note that both the number $\widehat{J}$ and the elements of the partition $\{\widehat{\mathcal{B}}_1,\ldots,\widehat{\mathcal{B}}_{\widehat{J}}\}$ are random quantities, and that, by construction, no two consecutive entries of $\widehat{\nu}$ are identical. Using (2.3), the individual entries of the vector $\widehat{\nu}$ can be obtained explicitly, as shown next.

LEMMA 2.1. *Let $\widehat{\nu} \in \mathbb{R}^{\widehat{J}}$ satisfy (2.3) and $\widehat{b}_j = |\widehat{\mathcal{B}}_j|$ for $1 \leq j \leq \widehat{J}$. Then,*

$$\widehat{\nu}_j = \frac{1}{\widehat{b}_j} \sum_{i \in \widehat{\mathcal{B}}_j} y_i + \widehat{c}_j,$$

*where*

(2.4) $$\widehat{c}_1 = \begin{cases} -\dfrac{\lambda_{2,n}}{\widehat{b}_j}, & \text{if } \widehat{\nu}_2 - \widehat{\nu}_1 > 0, \\ \dfrac{\lambda_{2,n}}{\widehat{b}_j}, & \text{if } \widehat{\nu}_2 - \widehat{\nu}_1 < 0, \end{cases}$$

(2.5) $$\widehat{c}_{\widehat{J}} = \begin{cases} \dfrac{\lambda_{2,n}}{\widehat{b}_j}, & \text{if } \widehat{\nu}_J - \widehat{\nu}_{J-1} > 0, \\ -\dfrac{\lambda_{2,n}}{\widehat{b}_j}, & \text{if } \widehat{\nu}_J - \widehat{\nu}_{J-1} < 0, \end{cases}$$

*and, for $1 < j < \widehat{J}$,*

(2.6) $$\widehat{c}_j = \begin{cases} \dfrac{2\lambda_{2,n}}{\widehat{b}_j}, & \text{if } \widehat{\nu}_{j+1} - \widehat{\nu}_j > 0, \widehat{\nu}_j - \widehat{\nu}_{j-1} < 0, \\ -\dfrac{2\lambda_{2,n}}{\widehat{b}_j}, & \text{if } \widehat{\nu}_{j+1} - \widehat{\nu}_j < 0, \widehat{\nu}_j - \widehat{\nu}_{j-1} > 0, \\ 0, & \text{if } (\widehat{\nu}_j - \widehat{\nu}_{j-1})(\widehat{\nu}_{j+1} - \widehat{\nu}_j) = 1. \end{cases}$$

By Proposition 1 in Friedman et al. (2007), the fused lasso estimator is obtained by soft-thresholding of the individual coordinates of $\widehat{\mu}^F$, so that we immediately obtain the next result.

COROLLARY 2.2. *The fused lasso estimator $\widehat{\mu}^{\mathrm{FL}}$ is*

(2.7) $$\widehat{\mu}_i^{\mathrm{FL}} = \begin{cases} \widehat{\mu}_i^F - \lambda_{1,n}, & \widehat{\mu}_i^F \geq \lambda_1, \\ 0, & |\widehat{\mu}_i^F| < \lambda_{1,n}, \\ \widehat{\mu}_i^F + \lambda_{1,n}, & \widehat{\mu}_i^F \leq -\lambda_1, \end{cases} \quad i = 1,\ldots,n,$$

*where $\widehat{\mu}^F$ is the fusion estimator.*

REMARKS.

1. As is apparent from Lemma 2.1, the individual blocks found by the fusion solution $\widehat{\mu}^F$ are each biased by a term whose magnitude depends directly on the regularization parameter $\lambda_{2,n}$ and, inversely, on the size of the estimated block itself. That is, the larger the estimated blocks the smaller the effect of the bias. This term is simply a vertical shift, which is positive if the block is a local maximum, negative if it is a local minimum, and is zero otherwise. See Figure 3. It is worth pointing out that, as expected, the solution obtained using the taut-string algorithm of Davies and Kovac (2001a) with global squeezing exhibits exactly the same behavior, with the magnitude of the vertical shift being controlled by the size of the tube around the integrated process instead of the penalty term $\lambda_{2,n}$.
2. The regularization parameter $\lambda_{1,n}$ modulates the magnitude of the sparsity penalty and induces some bias effect as well. However, unlike the bias determined by the total variation penalty, this second type of bias is of the same magnitude for all the nonzero coordinates, a fact that can be seen directly from (2.7). An easy fix, which is considered in Section 2.3, is to adaptively penalize the estimated blocks differently, depending on their sizes, with larger blocks penalized less.

2.2. *Sparsistency for the fused lasso.* In this section, we provide conditions under which the block partition $\{\mathcal{B}_1^0, \ldots, \mathcal{B}_{J_0}^0\}$ and the block sparsity pattern $\mathcal{JS}(\mu^0)$ of $\mu^0$ can be estimated consistently [see (1.1)] by the fused lasso procedure. We break down our analysis into two parts, dealing separately with the fusion estimator $\widehat{\mu}^F$ first, which can be used to recover $\{\mathcal{B}_1^0, \ldots, \mathcal{B}_{J_0}^0\}$, and then with the fused lasso solution $\widehat{\mu}^{FL}$, from which the set $\mathcal{JS}(\mu^0)$ can be estimated. In Section 2.3, we show how this second task can be accomplished more effectively by a modified version of the fused lasso estimator.

2.2.1. *Recovery of true blocks by fusion only.* We first derive sufficient conditions for the fusion estimator to recover correctly the block partition of $\mu^0$. Let $\mathcal{J}_0 = \mathcal{J}(\mu^0)$ be the set of jumps of $\mu^0$ and $J_0 = |\mathcal{J}_0| + 1$ the cardinality of the associated block partition. Similarly, let $\widehat{\mathcal{J}} = \mathcal{J}(\widehat{\mu}^F)$ be the set of jumps for the fusion estimate given in (2.3).

THEOREM 2.3. *Assume* (E) *and (1.1). If, for some $\delta > 0$:*

1. $\frac{\lambda_{2,n}}{\sigma_n} \to \infty$ *and* $\frac{\lambda_{2,n}}{\sigma_n \sqrt{\log(n-J_0)}} > \frac{1}{2\sqrt{2}}(1+\delta)$,
2. $\frac{b_{\min}^0 \alpha_n}{\sigma_n} \to \infty$, $\frac{b_{\min}^0 \alpha_n}{\sigma_n \sqrt{\log J_0}} > \sqrt{16}(1+\delta)$ *and* $\lambda_{2,n} < b_{\min}^0 \frac{\alpha_n}{4}$,

*where $\alpha_n = \min_{i \in \mathcal{J}_0} |\mu_i^0 - \mu_{i-1}^0|$ and $b_{\min}^0 = \min_{1 \le j \le J_0} b_j^0$. Then,*

$$(2.8) \quad \lim_n \mathbb{P}(\{\widehat{\mathcal{J}} = \mathcal{J}_0\} \cap \{\operatorname{sgn}(\widehat{\mu}_i^F - \widehat{\mu}_{i-1}^F) = \operatorname{sgn}(\mu_i^0 - \mu_{i-1}^0), \forall i \in \mathcal{J}_0\}) = 1.$$

**Fusion Estimates**

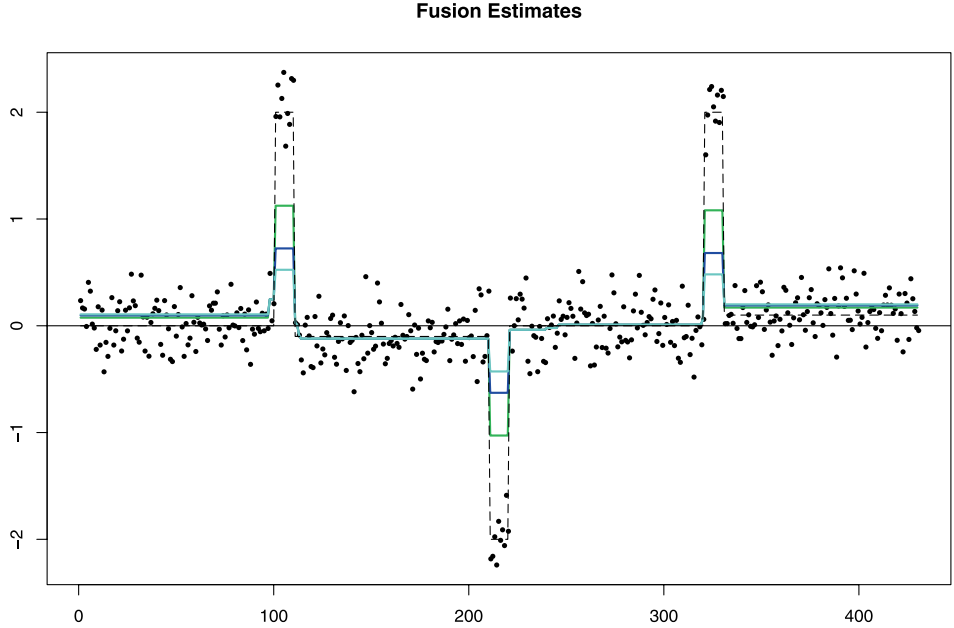

Fig. 3. *Different fusion estimates for the data described in Section 2.4. The dashed line corresponds to the true mean vector, while the three lines correspond to the fusion estimates with different regularization parameters.*

Remarks.

1. In the proof of Theorem 2.3, instead of Slepian's inequality, one could use Markov's inequality and well-known bounds on the supremum of centered sub-Gaussian vectors [see, e.g., Lemma 2.3 in Massart (2007)] to derive slightly stronger sufficient conditions for (2.8), which however hold for the larger class of sub-Gaussian errors. We give the following conditions without a proof:

   (a) $\lim_n \frac{\sigma_n \sqrt{2\log|\mathcal{J}(\mu_0)|} + 2\lambda_{2,n}}{b^0_{\min} \alpha_n} = 0$,

   (b) $\lim_n \frac{\lambda_{2,n}}{\sigma_n \sqrt{\log|\mathcal{J}_0^c|}} = \infty$.

   Furthermore, the errors need not be identically distributed. In fact, the proof of the theorem holds almost unchanged if, for example, one only assumes that the individual variances are of order $O(1/\sqrt{n})$.

2. Equation (2.8) actually implies not only that $\mathcal{J}_0$ can be consistently estimated, but also that the true signs of the jumps can be recovered with overwhelming probability, a feature known in the lasso literature as sign consistency [see, e.g., Wainwright (2006) and Zhao and Yu (2006)]. In the present settings, sign consistency of the fusion estimate implies the following desirable feature of $\widehat{\mu}^F$:

COROLLARY 2.4. *The fusion estimator $\widehat{\mu}^F$ can consistently recover the local maxima and local minima of $\mu^0$.*

3. The magnitude $\alpha_n$ of the smallest jumps of $\mu^0$ is a fundamental quantity, whose asymptotic behavior determines whether recovery of the true blocks obtains. In particular, if $\alpha_n$ vanishes at a rate faster than $\sqrt{b^0_{\min}}/\sigma_n$, then no recovery is possible. In a way, this guarantees some form of asymptotic distinguishability that prevents adjacent blocks from looking too similar.
4. The larger the minimal size of a block $b^0_{\min}$, the easier the recovery of the blocks by fusion.

2.2.2. *Recovery of true blocks and true nonzero coordinates by the fused lasso.* Let $\mathcal{JS}_0 = \mathcal{JS}(\mu^0)$ be set of nonzero blocks of $\mu^0$ and $K_0 = |\mathcal{JS}_0|$ its cardinality. Let $\widehat{\mathcal{JS}} = \mathcal{JS}(\widehat{\mu}^{\mathrm{FL}})$ be the equivalent quantity defined using the fused lasso estimate $\widehat{\mu}^{\mathrm{FL}}$. Consider the event

$$\mathcal{R}_{1,n} = \{\mathcal{JS}_0 = \widehat{\mathcal{JS}}\} \cap \{\operatorname{sgn}(\widehat{\nu}_j) = \operatorname{sgn}(\nu^0_j), \forall j \in \mathcal{JS}_0\}$$

that soft-thresholding $\widehat{\mu}^F$ with penalty $\lambda_{1,n}$ will return the nonzero blocks of $\mu^0$.

THEOREM 2.5. *If the conditions of Theorem 2.3 are satisfied and, for some $\delta > 0$:*

1. $\frac{\lambda_{1,n}\sqrt{b^0_{\min}}}{\sigma_n} \to \infty$ *and* $\frac{\lambda_{1,n}\sqrt{b^0_{\min}}}{\sigma_n\sqrt{\log(J_0 - K_0)}} > 2\sqrt{2}(1+\delta)$;
2. $\frac{2\lambda_{2,n}}{b^0_{\min}} < \frac{\lambda_{1,n}}{2}$, *for all $n$ large enough;*
3. $\frac{\rho_n\sqrt{b^0_{\min}}}{\sigma_n} \to \infty$, $\frac{\rho_n\sqrt{b^0_{\min}}}{\sigma_n\sqrt{\log K_0}} > \sqrt{18}(1+\delta)$ *and* $\lambda_{1,n} < \frac{\rho_n}{3}$ *for all $n$ large enough;*
4. $\frac{2\lambda_{2,n}}{b^0_{\min}} < \frac{\rho_n}{3}$, *for all $n$ large enough,*

*where $\rho_n = \min_{j \in \mathcal{K}_0} |\nu^0_j|$; then,*

$$\lim_n \mathbb{P}(\mathcal{R}_{1,n}) = 1.$$

REMARKS.

1. As was the case for Theorem 2.3, the assumption of Gaussian errors is not essential and can be relaxed, and, in fact, Remark 1 above still applies.
2. The previous result implies that the fused lasso is not only consistent but, in fact, sign consistent, so that the signs of the nonzero blocks are estimated correctly.

3. The magnitude $\rho_n$ of the smallest nonzero block value cannot decrease to zero too fast, otherwise the sparsity pattern cannot be fully recovered, just as we pointed out above in Remark 3 for the fusion solution.
4. The conditions of Theorem 2.5 appear to be quite cumbersome for two main reasons. First, the regularization parameters $\lambda_{1,n}$ and $\lambda_{2,n}$ interact with each other. As a result, it appears necessary to impose assumption 2 in order to guarantee that the two different bias terms they each determine will not disrupt the recovery process. Secondly, one has to keep track of the size $b_{\min}^0$ of the minimal block. This additional bookkeeping is due to the fact that the sparsity penalty is enforced globally, in the sense that all coordinates are penalized in equal amount, thus ignoring the fact that longer blocks require less regularization (see Remark 1 after Lemma 2.1).

2.3. *The fused adaptive lasso: Sparsistency and an oracle inequality.* Motivated by the stringent nature of the conditions of Theorem 2.5, below we propose a refinement of the fused lasso estimator, which we call the *fused adaptive lasso*. Overall, this slightly different estimator enjoys better asymptotic properties than the fused lasso, at no additional complexity cost.

The fused adaptive lasso is obtained with the following two-step procedure:

1. Fusion step. Compute the fusion solution $\widehat{\mu}^F$ using the fusion regularization parameter $\lambda_{2,n}$, as in (2.2), and the corresponding block-partition $(\widehat{\mathcal{B}}_1, \ldots, \widehat{\mathcal{B}}_{\widehat{J}})$ [see (2.3)]. Obtain

$$\widehat{\mu}^{\mathrm{AF}} = \sum_{j=1}^{\widehat{J}} \bar{y}_j 1_{\widehat{\mathcal{B}}_j}, \tag{2.9}$$

where

$$\bar{y}_j = \frac{1}{\widehat{b}_j} \sum_{i \in \widehat{\mathcal{B}}_j} y_i, \qquad 1 \le j \le \widehat{J}.$$

2. Adaptive lasso step. Compute the fused adaptive lasso solution

$$\widehat{\mu}^{\mathrm{FAL}} = \operatorname*{arg\,min}_{\mu \in \mathbb{R}^n} \|\widehat{\mu}^{\mathrm{AF}} - \mu\|_2^2 + \sum_{i=1}^n \tilde{\lambda}_i |\mu_i|, \tag{2.10}$$

where the $n$-dimensional random vector $\tilde{\lambda}$ of penalties is

$$\lambda = \lambda_1 \sum_{j=1}^{\widehat{J}} \frac{1}{\sqrt{\widehat{b}_j}} 1_{\widehat{\mathcal{B}}_j} \tag{2.11}$$

with $\lambda_{1,n}$ as the $\ell_1$ regularization parameter.

REMARKS.

1. The fused adaptive lasso differs from the fused lasso in two fundamental aspects. First, as easily seen from (2.9), the bias term in the fusion solution due to the terms $c_j$, which depends on the regularization parameter $\lambda_{2,n}$, is absent (see Lemma 2.1). Equivalently, the fusion estimator is only used to estimate the block partition of $\mu^0$, and, provided this estimate is correct, the block values are estimated unbiasedly with the sample averages. Using the fusion procedure as an estimator of the block partition has the other advantage of decoupling the estimation from the model selection problem, thus freeing, to some extent, the user from the task of carefully choosing an optimal penalty $\lambda_{2,n}$. In fact, recovery of the true partition can be obtained even if the problem is overpenalized, and, therefore, the resulting estimator $\widehat{\mu}^F$ is highly biased.

   Secondly, the penalty terms used for thresholding individual blocks are rescaled by the squared root of the length of the estimated blocks. The rationale for using this rescaling is very simple. In fact, suppose that, for some $j_1, j_2$, $b_{j_1} \gg b_{j_2}$. Since the variance of the $j$th block average $\bar{y}_j$ is $\frac{\sigma_n^2}{b_j}$, $\bar{y}_{j_1}$ has a much smaller standard error than $\bar{y}_{j_2}$ and, therefore, should be penalized less heavily. The adequate reduction in the sparsity penalty of $\bar{y}_{j_1}$ versus $\bar{y}_{j_2}$ is precisely the difference in their standard errors, hence the choice of rescaling by the square root of the block lengths. The advantage of adaptively thresholding the block values in this manner is that the procedure will be more effective at identifying longer nonzero blocks whose values are quite close to 0.

   In Section 2.4 we explain both these improvements concretely with a numerical example.

2. In step 2 the vector $\widehat{\mu}$ is straightforward to compute via soft-thresholding of the individual coordinates of $\widehat{\mu}^{\mathrm{AF}}$ with coordinate-dependent thresholds

$$\widehat{\mu}_i^{\mathrm{FAL}} = \begin{cases} \widehat{\mu}_i^{\mathrm{AF}} - \lambda_i, & \widehat{\mu}_i^{\mathrm{AF}} \geq \lambda_i, \\ 0, & |\widehat{\mu}_i^{\mathrm{AF}}| < \lambda_i, \\ \widehat{\mu}_i^{\mathrm{AF}} + \lambda_i, & \widehat{\mu}_i^{\mathrm{AF}} \leq -\lambda_i, \end{cases} \qquad 1 \leq i \leq n.$$

3. Instead of the soft-thresholded block estimate of step 2, one may consider instead the corresponding estimate $\widetilde{\mu}$ based on the hard-threshold where

$$\widetilde{\mu}_i = \widehat{\mu}_i^{\mathrm{AF}} 1\{|\widehat{\mu}_i^{\mathrm{AF}}| \geq \lambda_i\}, \qquad 1 \leq i \leq n.$$

One of the asymptotic advantages of the fused adaptive lasso versus the ordinary fused lasso is that block recovery obtains under milder conditions than Theorem 2.5, without the need to consider the fusion penalty parameter $\lambda_{2,n}$ and the length of the minimal block. In some sense, the fused adaptive lasso can adapt more flexibly to the block sparsity than the fused lasso.

PROPOSITION 2.6. *Assume that the conditions of Theorem 2.3 are satisfied. Then,*

$$\lim_n \mathbb{P}\{\mathcal{R}_{1,n}\} = 1,$$

*if, for some $\delta > 0$,*

1. $\frac{\lambda_{1,n}}{\sigma_n} \to \infty$ *and* $\frac{\lambda_{1,n}}{\sigma_n \sqrt{\log(J_0 - K_0)}} > \sqrt{2}(1+\delta)$;
2. $\frac{\rho_n}{\sigma_n} \to \infty$, $\frac{\rho_n}{\sigma_n \sqrt{\log K_0}} > 2\sqrt{2}(1+\delta)$ *and* $\lambda_{1,n} < \frac{\rho_n}{2}$ *for all $n$ large enough,*

*where $\rho_n = \min_{j \in \mathcal{K}_0} |\nu_j^0|$.*

A second advantage of the fused adaptive lasso stems from the oracle property derived below. Consider the ideal situation where we have access to an oracle who lets us know the $K^0$ sets $\mathcal{B}_{j_k}^0$, $k = 1, \ldots, K^0$, of the true block partition of $\mu^0$ for which $|\nu_{j_k}^0| > \sigma_n / \sqrt{b_{j_k}^0}$. Notice that, from this information, one can recover the true partition. The oracle estimate $\widehat{\mu}^O$ is the vector with coordinates

$$\widehat{\mu}_i^O = \begin{cases} \frac{1}{b_{j_k}^0} \sum_{z \in \mathcal{B}_{j_k}^0} y_z, & \text{if } i \in \mathcal{B}_{j_k}, \\ 0, & \text{otherwise.} \end{cases}$$

This procedure amounts to setting to 0 the estimates for the coordinates belonging to the blocks whose true mean value is smaller than $\sigma_n / \sqrt{b_j^0}$. The corresponding ideal risk is

$$\begin{aligned}
\mathbb{E}\|\widehat{\mu}^O - \mu^0\|_2^2 &= \sum_i \sum_{j_k} 1\{i \in \mathcal{B}_{j_k}^0\} \min\left\{\frac{\sigma_n^2}{b_{j_k}^0}, (\nu_j^0)^2\right\} \\
&= K_0 \sigma_n^2 + \sum_{j \notin \mathcal{JS}_0} b_j^0 (\nu_j^0)^2.
\end{aligned}$$
(2.12)

Note, in particular, that

$$\mathbb{E}\|\widehat{\mu}^O - \mu^0\|_2^2 \leq \sum_i \min\{\sigma_n^2, \mu_i^2\}$$

with equality if and only if $b_j^0 = 1$ for all $j$, where the expression on the right-hand side is the ideal risk for the oracle estimator based on thresholding of individual coordinates rather than of blocks. Therefore, if $\mu^0$ has a block structure, as is assumed here, this different oracle will be able to achieve a smaller ideal risk.

Before stating our oracle result, we need some additional notation. Recall that any $\mu \in \mathbb{R}^n$ can always be written as

$$\mu = \sum_{j=1}^{J} \nu_j 1_{\mathcal{B}_j} \tag{2.13}$$

for some (possibly trivial) block partition $(\mathcal{B}_1, \ldots, \mathcal{B}_J)$ of $\{1, \ldots, n\}$, with $J \leq n$. Let $\mu^1$ and $\mu^2$ be vectors in $\mathbb{R}^n$ with block partitions $\{\mathcal{B}_1^1, \ldots, \mathcal{B}_{J_1}^1\}$ and $\{\mathcal{B}_1^2, \ldots, \mathcal{B}_{J_2}^2\}$, respectively, where $J_1, J_2 \leq n$. Then, they satisfy (2.13), for some vectors $\nu^1 \in \mathbb{R}^{J_1}$ and $\nu^2 \in \mathbb{R}^{J_2}$, respectively. Let $\{\mathcal{L}_1, \ldots, \mathcal{L}_m\}$ be the partition of $\{1, \ldots, n\}$ obtained as the refinement of the block partitions of $\mu^1$ and $\mu^2$, that is, for every $l = 1, \ldots, m$, $\mathcal{L}_l = \mathcal{B}_{j_1}^1 \cap \mathcal{B}_{j_2}^2$, for some $j_1$ and $j_2$. We define the quantity

$$\mathcal{JS}(\mu^1; \mu^2) = \{l : \mathcal{L}_l = \mathcal{B}_{j_1}^1 \cap \mathcal{B}_{j_2}^2, \nu_{j_1}^1 \neq 0\}.$$

THEOREM 2.7. *Assume that $\mu^0$ satisfies (1.1) and that*

$$\alpha_n = O\left(\sqrt{\frac{\log n}{n}}\right). \tag{2.14}$$

*Let $\sigma_n^2 = \frac{\sigma^2}{n}$, $\lambda_{2,n} = A\sqrt{\sigma_n^2 \log n}$, with $A > 0$ such that $\lambda_{2,n}\alpha_n < 1/4$ and $\lambda_{1,n} = 2\sqrt{\sigma_n^2 \log \widehat{J}}$, where $\widehat{J}$ is obtained by solving the fusion problem (2.2) in the first step of the adaptive fused-lasso procedure. For any vector $\mu \in \mathbb{R}^n$, set*

$$V(\mu) = 32|\mathcal{JS}(\mu; \mu^0)|\sigma_n^2 \log J_0.$$

*Then, for any $\delta \in [0, 1)$,*

$$\lim_n \mathbb{P}\left\{\|\widehat{\mu}^{\mathrm{FAL}} - \mu_0\|_2^2 \leq \frac{2+\delta}{2-\delta} \inf_{\mu \in \mathbb{R}^n} \{V(\mu) + \|\mu - \mu^0\|_2^2\}\right\} = 1. \tag{2.15}$$

REMARKS.

1. The assumption in (2.14) stems from Theorem 2.3 and is crucial in our proof, as it guarantees that recovery of the true block partition of $\mu^0$ by fusion, which is necessary for mimicking the oracle solution $\widehat{\mu}^O$, is feasible. It essentially allows for consecutive blocks to differ by a vanishing quantity of smaller order than $\sqrt{\log n/n}$. If the minimal jump size is bounded away from zero, uniformly in $n$, then the condition $\lambda_{2,n}\alpha_n < 1/4$ is redundant.

2. The proof of Theorem 2.7 shows that $V(\mu)$ is minimized by vectors such that

$$|\mathcal{JS}(\mu;\mu^0)| = |\mathcal{JS}(\mu^0)| = K^0;$$

that is, vectors whose block partition matches the the true block partition. Therefore, (2.15) shows that the adaptive fused-lasso achieves the same oracle rates granted by ideal risk (2.12) up to a term that is logarithmic in $J_0$.
3. If it is further assumed that $\|\mu^0\|_\infty < C$ uniformly in $n$, for some constant $C$, the result (2.15) can be strengthened to

$$\mathbb{E}\|\widehat{\mu} - \mu_0\|_2^2 \leq \frac{2+\delta}{2-\delta} \inf_{\mu \in \mathbb{R}^n} \{V(\mu) + \|\mu - \mu^0\|_2^2\} + o(1).$$

2.4. *A toy example.* We discuss a stylized numerical example for the purpose of clearly illustrating the two advantages of the fused adaptive lasso, namely the use of the fusion penalty only for recovering the true block partition and the block-dependent rescaling of the lasso penalty. See Remark 1 before Proposition 2.6 for details.

We simulate one sample according to the model

$$y_i = \mu_i^0 + \varepsilon_i,$$

where

$$\mu_i^0 = \begin{cases} 0, & 1 \leq i \leq 100, \\ 2, & 101 \leq i \leq 110, \\ -0.1, & 111 \leq i \leq 210, \\ -2, & 211 \leq i \leq 220, \\ 0, & 221 \leq i \leq 320, \\ 2, & 321 \leq i \leq 330, \\ 0.1, & 331 \leq i \leq 430, \end{cases}$$

and the errors are independent Gaussian variables with mean zero and standard deviation $\sigma = 0.2$. Figure 1 shows the data along with the true signal. Notice that some of the coordinates of $\mu^0$ are in absolute value less than $\sigma$, a fact that, as we will see, if $\mu^0$ were not blocky, would make the recovery of those coordinates infeasible. Figure 3 portrays the simulated data and three fusion estimates $\widehat{\mu}^F$, each of them solving (2.2) for three different values of $\lambda_{2,n}$: 4.8, 6.8 and 7.8. The dashed line corresponds to the true mean vector $\mu^0$. The excessive amount of penalization is apparent from the large bias in all these estimates, especially in the smaller blocks. Nonetheless, the block partitions that each of these estimates produce match, in fact, very closely the true block partition.

Figure 4 shows the modified fusion estimate $\widehat{\mu}^{\mathrm{AF}}$ given in (2.9) using the fusion estimate from Figure 3 with the largest amount of bias, along with

**Modified Fusion Estimate**

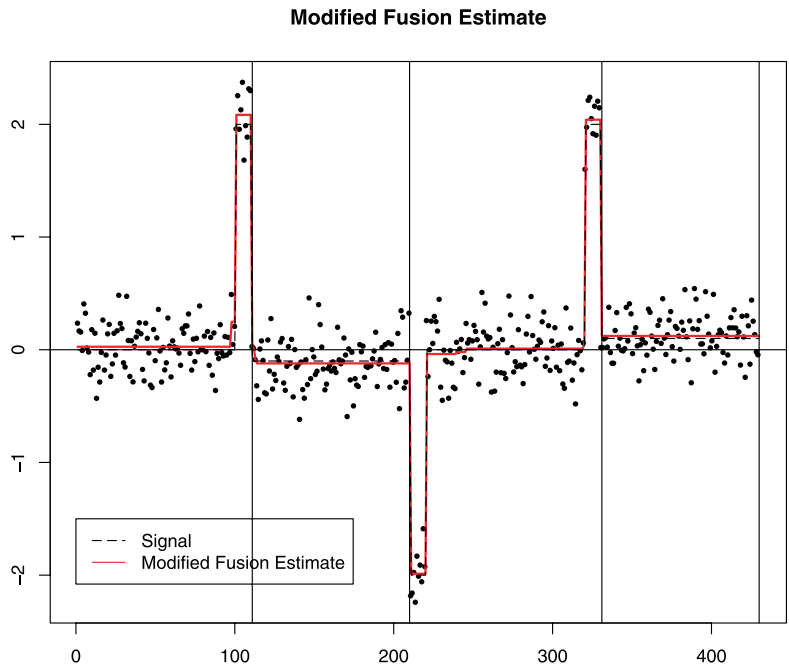

Fig. 4. *The modified fusion estimate $\widehat{\mu}^{\mathrm{AF}}$ of (2.9), using the fusion estimate from Figure 3 with the lowest total variation. The dashed gray line, which is almost indistinguishable from the estimate, is the true signal $\mu^0$. The vertical lines enclose the third and seventh blocks, whose value is in magnitude half the standard deviation of the errors.*

the true mean vector $\mu^0$, displayed as a dashed line. Because the block partition was estimated correctly, the estimate $\widehat{\mu}^{\mathrm{AF}}$ is almost indistinguishable from the true vector $\mu^0$. For this particular dataset, the adaptive lasso step would set to zero correctly the first and fifth block, but not the third and seventh blocks, which in Figure 4 are enclosed by black vertical lines. In fact, although the true value of those blocks is in magnitude half the standard deviation of the errors, $\sigma$, the standard error for both the block estimates is roughly $\sigma/10$. This is taken into account in the adaptive lasso step, but not in the lasso step, where even the ideal soft threshold, that is $\sigma$, would be too high, thus incorrectly setting to zero both of these blocks.

Finally, we simulated 1000 datasets according to the model described here and computed the empirical mean squared errors for the fused adaptive lasso estimates, using for the penalty terms the values indicated in Theorem 2.7. Figure 5 shows the histogram of the empirical mean squared errors, with the vertical line representing the true mean squared error $\frac{1}{n}\mathbb{E}\|y - \mu^0\|^2$, namely $\sigma^2$. Notice how the empirical mean squared errors are larger then the true value, the usual price paid for adaptivity.

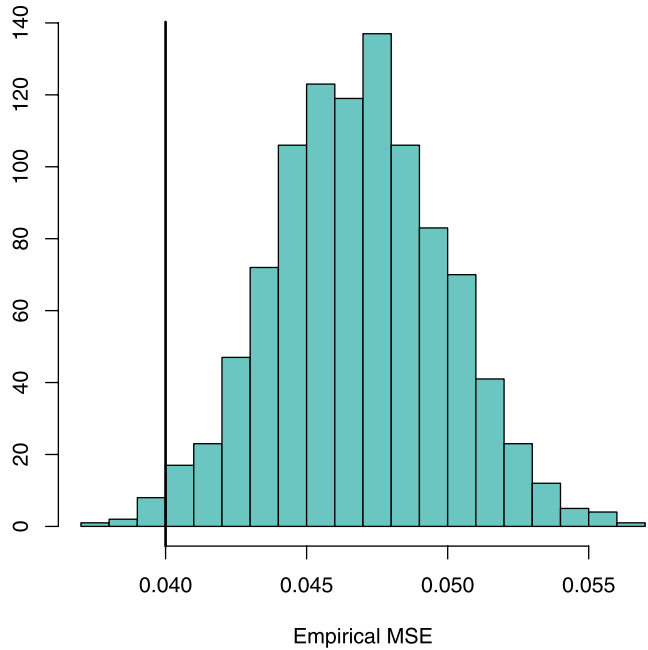

Fig. 5. *Distributions of the empirical mean squared errors from 1000 simulations from the model described on Section 2.4 using the fused adaptive lasso with penalty parameters chosen according to Theorem 2.7. The vertical line represents $\sigma^2$.*

2.5. *How to choose $\lambda_{1,n}$ and $\lambda_{2,n}$.* From the practical standpoint, the choice of the regularization parameters is crucial. For the fused adaptive lasso, one can infer from the proof of Theorem 2.7 that the optimal choice for the vector of lasso penalty terms $\lambda$ is given by

$$(2\sigma_n\sqrt{\log \widehat{J}}) \sum_{j=1}^{\widehat{J}} \frac{1}{\sqrt{\widehat{b}_j}} 1_{\widehat{\mathcal{B}}_j},$$

with $1_{\widehat{\mathcal{B}}_j}$ denoting the indicator vector of the estimated block $\widehat{\mathcal{B}}_j$, $1 \leq j \leq \widehat{J}$. This choice corresponds to soft-thresholding $\widehat{J}$ independent Gaussian variables with variances $\frac{\sigma_n^2}{\sqrt{\widehat{b}_j}}$, $j = 1, \ldots, \widehat{J}$.

Admittedly, for the total variation penalty term $\lambda_{2,n}$ the theoretical results presented here, being of asymptotic nature, may not directly lead to procedures that are effective in practice, unless $n$ is very large. Choosing optimal values for the penalty parameters remains an important open problem in much of the penalized least-squares literature, where the theoretical (e.g., asymptotic) results may offer little guidance in practice. Cross validation is certainly a viable way of choosing both $\lambda_{1,n}$ and $\lambda_{2,n}$, as recommended in Friedman et al. (2007), and as is almost exclusively done in practice (although it remains to be seen whether this procedure leads to optimal estimators). Nonetheless, an automatic procedure for choosing $\lambda_{2,n}$ that exhibits reasonable performance still eludes us. However, our theoretical analysis, and the toy example presented above, shows that, if the signal is comprised mostly of long blocks, a large value of $\lambda_{2,n}$ will lead to accurate estimates of the block partition, and the results should be relatively robust to different choices.

An interesting possibility suggested by a referee, which is beyond the scope of this article, is to replace the overall total variation parameter $\lambda_{2,n}$ with a series of data-driven parameters, one for each term of the total variation seminorm. Specifically, one can consider the penalized problem

$$(2.16) \qquad \arg\min_{\mu \in \mathbb{R}^n} \left\{ \sum_{i=1}^n (y_i - \mu_i)^2 + 2 \sum_{i=2}^n \lambda_{2,i} |\mu_i - \mu_{i-1}| \right\},$$

where $\{\lambda_{2,i}, i = 2, \ldots, n\}$ are possibly different coefficients that modulate the effect of the total variation penalty at different locations along the signal, so that the solution is more robust to spurious local extreme due to unusually large errors. In fact, as pointed out by Davies and Kovac (2001b), the taut-string algorithm with local squeezing approximates the solution to this problem. Although local squeezing increases the complexity of the algorithm, it has been shown to enjoy a better performance than the problem with an omnibus total variation penalty. The choice of the regularization parameters $\{\lambda_{2,i}, i = 2, \ldots, n\}$ can be done iteratively, starting with all $\lambda_{2,i}$'s being identical and very large (thus producing an estimate with constant entries) and then, at every step, shrinking them differently based on the features of the residuals, such as the multiresolution coefficients as defined in Davies and Kovac (2001a).

**3. Sieve methods.** In this section, we study the rates of convergence for the sieve least squares solutions (1.3) and (1.4). For convenience, consistency is measured with respect to the normalized Euclidean norm $\|x\|_n = \frac{1}{\sqrt{n}} \sqrt{\sum_{i=1}^n x_i^2}$. Accordingly, we change our assumption on the errors as follows:

(E′) The errors $(\varepsilon_1, \ldots, \varepsilon_n)$ are independent sub-Gaussian variables with variances bounded by $\sigma^2$, uniformly in $n$.

Notice that the results and settings of previous sections can be adapted in a straightforward way to the present framework.

We first study the estimator given in (1.3). To that end, consider the class of vectors

$$\mathcal{C}_{\text{TV}}(T_n) = \{\mu \in \mathbb{R}^n : \|\mu\|_{\text{TV}} \leq T_n, \|\mu\|_\infty \leq C\},$$

where $C$ is a finite constant that does not depend on $n$, and the $\ell_1$-ball of radius $L_n$

$$\mathcal{C}_{\ell_1}(L_n) = \{\mu \in \mathbb{R}^n : \|\mu\|_1 \leq L_n\}$$

with both numbers $T_n$ and $L_n$ being allowed to grow unboundedly with $n$. Then, we can rewrite (1.3) as

$$\widehat{\mu}^{TL} = \underset{\mu \in \mathcal{C}_{\text{TV}}(T_n) \cap \mathcal{C}_{\ell_1}(L_n)}{\arg\min} \|y - \mu\|_2^2.$$

Below, we derive the consistency rate for $\widehat{\mu}^{TL}$ in terms of the sequences $T_n$ and $L_n$ by dealing separately with the two sieves.

THEOREM 3.1. *Assume* (E') *and* $\mu_n^0 \in \mathcal{C}_{\ell_1}(L_n) \cap \mathcal{C}_{\text{TV}}(T_n)$. *Let*

(3.1) $$\widehat{\mu}^T = \inf_{\mu \in \mathcal{C}_{\text{TV}}(T_n)} \|y - \mu\|_2^2$$

*and*

$$\widehat{\mu}^L = \underset{\mu \in \mathcal{C}_{\ell_1}(L_n)}{\arg\min} \|y - \mu\|_2.$$

*Then,*

$$\|\widehat{\mu}^T - \mu^0\|_n = O_P(T_n^{1/3} n^{-1/3}),$$

*so that $\widehat{\mu}^T$ is consistent provided that $T_n = o(n)$, and*

(3.2) $$\|\widehat{\mu}^L - \mu^0\|_n = O_P\left(\sqrt{\frac{L_n (\log n)^{3/2}}{n}}\right),$$

*so that $\widehat{\mu}^L$ is consistent provided that*

$$L_n = o\left(\frac{n}{(\log n)^{3/2}}\right).$$

*As a result,*

(3.3) $$\|\widehat{\mu}^{TL} - \mu^0\|_n = O_P\left(\frac{L_n (\log n)^{3/2}}{n} \wedge \left(\frac{T_n}{n}\right)^{1/3}\right).$$

REMARKS.

1. It appears that the requirement for the vectors in $\mathcal{C}_{\text{TV}}(T_n)$ to be uniformly bounded cannot be relaxed without negatively affecting the rate of consistency or without introducing additional assumptions [see, e.g., Theorem 9.2 in van de Geer (2000)].

2. The rate of consistency for $\widehat{\mu}^F$ should be compared with the analogous rate derived in Theorem 9 of Mammen and van de Geer (1997) for the penalized version of the least squares problem (3.1).
3. The rate given in (3.2) is not the sharpest possible. In fact, an application of Theorem 5 of Donoho and Johnstone (1994) yields for $\widehat{\mu}^L$ the improved minimax rate

$$\sqrt{\frac{L_n}{n}}(\log n)^{1/4}$$

for the case of i.i.d. Gaussian errors, from which we can infer a maximal rate of growth $L_n = o(\frac{n}{\sqrt{\log n}})$.
4. We make no claims that the rate given in equation (3.3), which is just the minimum of the rates for two separate sieve least squares problems, is sharp. Better rates may be obtained from better estimates of the metric entropy of the set $\mathcal{C}_{\ell_1}(L_n) \cap \mathcal{C}_{\mathrm{TV}}(T_n)$.
5. *On the relationship between $L_n$ and $T_n$.* The total variation and $\ell_1$ constraints are not independent of each other. One can easily verify that

$$T_n^{\max} \equiv \max_{x \in \mathcal{C}_{\ell_1}(L_n)} \|x\|_{\mathrm{TV}} = 2L_n.$$

On the other hand, every vector $x \in \mathbb{R}^n$ such that $\|x\|_{\mathrm{TV}} = T_n$ can be written as

$$x = m + t,$$

where $\|t\|_{\mathrm{TV}} = T_n$, $m = 1_n \bar{x}_n$, with $\bar{x}_n = \frac{1}{n}\sum_i x_i$, and $\frac{1}{n}\sum_i m_i t_i = 0$. Notice that $m$ can be estimated at the rate $\frac{1}{\sqrt{n}}$, so the convergence rates for $\widehat{\mu}^T$ depend on how well $t$ can be estimated. Next, notice that

$$L_n^{\max} \equiv \max_{x \in \mathcal{C}_{\mathrm{TV}}(T_n), x = m+t} \|t\|_1 = \frac{T_n}{2} \frac{n}{n-1},$$

where $m + t$ is the decomposition of $x$ discussed above. Therefore, over the set $\mathcal{C}_{\mathrm{TV}}(T_n) \cap \mathcal{C}_{\ell_1}(L_n)$, we obtain the relationship

(3.4) $$T_n^{\max} \sim 2L_n^{\max}.$$

Our final result concerns the estimator resulting from the nonconvex sieve least squares problem (1.4). Define the set

$$\mathcal{C}(S_n, J_n) = \{\mu \in \mathbb{R}^n : |\mathcal{S}_n(\mu)| \leq S_n\} \cap \{\mu \in \mathbb{R}^n : |\mathcal{J}_n(\mu)| + 1 \leq J_n\},$$

consisting of vectors in $\mathbb{R}^n$ that have at most $S_n$ nonzero coordinates and take on at most $J_n$ different values. For convenience, we further impose the following, fairly weak assumption, which does not preclude the coordinates of $\mu^0$ from becoming increasingly large in magnitude with $n$:

(R) the set $\mathcal{C}(S_n, J_n)$ is contained in a $S_n$-dimensional cube centered at the origin with volume $R_n$ such that

$$\log R_n = o(n).$$

THEOREM 3.2. *Assume* (E′) *and* (R) *and let* $\widehat{\mu}^{SJ} = \arg\min_{\mu \in \mathcal{C}(S_n, J_n)} \|y - \mu\|_2^2$.

1. *If* $S_n = o(\frac{n}{\log n})$, *then*

(3.5) $$\|\widehat{\mu}^{SJ} - \mu^0\|_n = O_P\left(\sqrt{\frac{J_n}{n}}\right).$$

2. *When* $S_n = n$, *(3.5) still holds, provided* $J_n = o(\frac{n}{\log n})$.

REMARKS.

1. The rate on $S_n$ is in accordance with the persistence rate derived in Greenshtein (2006), Theorem 1, for related least squares regression problems on sieves.
2. If $J_0$ is bounded, uniformly in $n$, the consistency rate we obtain is parametric. See Boysen et al. (2009) for a similar result.

**4. Discussion and future directions.** In this work, we tackle the task of estimating a blocky and sparse signal using three different methodologies, whose asymptotic properties we investigate. We study the fused lasso estimator proposed in Friedman et al. (2007) and a simple variant of it, with better properties. For both procedures, we provide conditions under which they recover with overwhelming probability as $n$ gets larger the block partition. We also study consistency rates of sieve least square problems under two types of constraints, one on the maximal radiuses of the $\ell_1$- and $\|\cdot\|_{\mathrm{TV}}$-balls, and the other on the maximal number of blocks and nonzero coordinates. Overall, these results complement each other in providing different types of asymptotic information for the task at hand and complement other analyses already existing in the statistical literature.

There are a number of generalizations of the results presented. We mention only the ones that seem the most natural to us. A first extension involves considering a corrupted version of a signal $\mu^0 \in \mathbb{R}^n \times \mathbb{R}^n$, corresponding to the problem of denoising a sparse, blocky image over a $n \times n$ grid, for which total variation methods have proven quite effective. Another interesting direction would be to assume a known slowly-varying variance function, for example, with given Lipschitz constant, and incorporate this information directly into the penalty functions in both the fusion and adaptive lasso steps. Furthermore, under this heteroschedastic scenario, one could first build a consistent estimator of the variance function and then, in the fusion step, use it to penalize the individual jumps adaptively. We think that our techniques and results can be directly generalized to study these more complex settings. Finally, we believe it would be quite valuable to investigate the possibility of building confidence balls and, in particular, confidence bands for the entire signal or for some of its local maxima or minima based on the estimators considered here.

## APPENDIX: PROOFS

LEMMA A.1. *Let $\|\cdot\|_{\mathrm{TV}}:\mathbb{R}^k \to \mathbb{R}$ be the fused penalty $\|x\|_{\mathrm{TV}} = \sum_{i=2}^k |x_i - x_{i-1}|$. Then, $\|\cdot\|_{\mathrm{TV}}$ is convex and, for any $x \in \mathbb{R}^k$, the subdifferential $\partial \|x\|_{\mathrm{TV}}$ is the set of all vectors $s \in \mathbb{R}^k$ such that*

$$
(A.1) \qquad s_i = \begin{cases} -w_2, & \text{if } i = 1, \\ w_i - w_{i+1}, & \text{if } 1 < i < k, \\ w_k, & \text{if } i = k, \end{cases}
$$

*where $w_i = \mathrm{sgn}(x_i - x_{i-1})$, for $2 \leq i \leq k$.*

PROOF. Let $L$ be a $(k-1) \times k$ matrix with entries $L_{i,i} = -1$ and $L_{i,i+1} = 1$ for $1 \leq i \leq (k-1)$ and 0 otherwise. Then, for any $x \in \mathbb{R}^k$, $\|x\|_{\mathrm{TV}} = \|Lx\|_1$. Convexity of $\|\cdot\|_{\mathrm{TV}}$ follows from the fact that it is the composition of a linear functional by the $\ell_1$ norm, which is convex. Next, by the definition of the subdifferential of the $\ell_1$ norm, for any $y \in \mathbb{R}^k$,

$$
(A.2) \qquad \|Ly\|_1 \geq \|Lx\|_1 + \langle L(y-x), w \rangle
$$

holds if and only if $w \in \mathcal{W}_x \subset \mathbb{R}^{k-1}$, where $\mathcal{W}_x$ is the set of all vectors $w$ such that $w_i = \mathrm{sgn}((Lx)_i)$. Equation (A.2) is equivalent to

$$
\|y\|_{\mathrm{TV}} \geq \|x\|_{\mathrm{TV}} + \langle y-x, s \rangle
$$

for each $k$-dimensional vector $s$ such that $s = L^\top w$ for some $w \in \mathcal{W}_x$. This set is described by (A.1) and is, therefore, $\partial \|x\|_{\mathrm{TV}}$. □

PROOF OF LEMMA 2.1. From the subgradient condition (2.1) with $\lambda_{1,n} = 0$, we obtain

$$
\widehat{\nu}_j = \frac{1}{\widehat{b}_j} \sum_{i \in \widehat{\mathcal{B}}_j} \widehat{\mu}_i^F = \frac{1}{\widehat{b}_j} \sum_{i \in \widehat{\mathcal{B}}_j} y_i - \frac{\lambda_{2,n}}{\widehat{b}_j} \sum_{i \in \widehat{\mathcal{B}}_j} s_i.
$$

Using (A.1), a simple telescoping argument leads to

$$
\sum_{i \in \mathcal{B}_j} s_i = w_{i_j} - w_{i_{j+1}} = \begin{cases} 2, & \text{if } (\widehat{\nu}_{j+1} - \widehat{\nu}_j) > 0, (\widehat{\nu}_j - \widehat{\nu}_{j-1}) < 0, \\ -2, & \text{if } (\widehat{\nu}_{j+1} - \widehat{\nu}_j) < 0, (\widehat{\nu}_j - \widehat{\nu}_{j-1}) > 0, \\ 0, & \text{if } (\widehat{\nu}_j - \widehat{\nu}_{j-1})(\widehat{\nu}_{j+1} - \widehat{\nu}_j) = 1, \end{cases}
$$

where $i_j = \min\{i : i \in \widehat{\mathcal{B}}_j\}$. This gives (2.6). It remains to consider the cases $j = 1$ and $j = \widehat{J}$. If $j = 1$, $\sum_{i \in \mathcal{B}_1} s_i = -w_{i_2}$, and if $j = \widehat{J}$, $\sum_{i \in \mathcal{B}_J} s_i = w_{i_J}$, form which (2.4) and (2.5) follow, respectively. □

PROOF OF THEOREM 2.3. Let

$$\text{(A.3)} \quad \mathcal{R}_{\lambda_{2,n}} = \{\widehat{\mathcal{J}} = \mathcal{J}_0\} \cap \{\text{sgn}(\widehat{\mu}_i^F - \widehat{\mu}_{i-1}^F) = \text{sgn}(\mu_i^0 - \mu_{i-1}^0), \forall i \in \mathcal{J}_0\}$$

and, for $2 \leq i \leq n$, let $d_i^0 = \mu_i^0 - \mu_{i-1}^0$, $\widehat{d}_i = \widehat{\mu}_i^F - \widehat{\mu}_{i-1}^F$ and $d_i^\varepsilon = \varepsilon_i - \varepsilon_{i-1}$. Using the subgradient conditions (A.1), the event $\mathcal{R}_{\lambda_{2,n}}$ occurs if and only if

$$d_i^\varepsilon = \lambda_{2,n}(2\,\text{sgn}(d_i^0) - \text{sgn}(\widehat{d}_{i-1}) - \text{sgn}(\widehat{d}_{i+1})) \qquad \forall i \notin \mathcal{J}_0,$$

where, for $x = 0$, $\text{sgn}(x)$ is the set $[-1, 1]$, and

$$|\widehat{d}_i| > 0 \qquad \forall i \in \mathcal{J}_0.$$

Next, in virtue of Lemma 2.1, on $\mathcal{R}_{\lambda_{2,n}}$ we can write

$$\widehat{d}_i = \frac{1}{b_{j(i)}^0} \sum_{k \in \mathcal{B}_{j(i)}^0} y_k + c_{j(i)}^0 - \frac{1}{b_{j(i-1)}^0} \sum_{k \in \mathcal{B}_{j(i-1)}^0} y_k - c_{j(i-1)}^0$$

$$= d_i^0 + \frac{1}{b_{j(i)}^0} \sum_{k \in \mathcal{B}_{j(i)}^0} \varepsilon_k - \frac{1}{b_{j(i-1)}^0} \sum_{k \in \mathcal{B}_{j(i-1)}^0} \varepsilon_k + c_{j(i)}^0 - c_{j(i-1)}^0,$$

where the index $j(i)$ identifies the block to which $i$ belongs; that is, $\mathcal{B}_{j(i)}^0$ is the block such that $i \in \mathcal{B}_{j(i)}$ for all $i = 1, \ldots, n$. Accordingly, $b_{j(i)} = |\mathcal{B}_{j(i)}^0|$ and $c_{j(i)}^0$ denotes the bias term in the fusion estimate as given in Lemma 2.1, with $\widehat{b}_j$ and $\widehat{\nu}_j$ replaced by $b_j^0$ and $\nu_j^0$, respectively, for $j = 1, \ldots, J_0$.

As a result, the event $\mathcal{R}_{\lambda_{2,n}}$ occurs in probability if both

$$\text{(A.4)} \quad \max_{i \notin \mathcal{J}_0} |d_i^\varepsilon| < \lambda_{2,n} |2\,\text{sgn}(d_i^0) - \text{sgn}(\widehat{d}_{i-1}) - \text{sgn}(\widehat{d}_{i+1})| < 4\lambda_{2,n}$$

and

$$\text{(A.5)} \quad \min_{i \in \mathcal{J}_0} \left| d_i^0 + \frac{1}{b_{j(i)}^0} \sum_{k \in \mathcal{B}_{j(i)}^0} \varepsilon_k - \frac{1}{b_{j(i-1)}^0} \sum_{k \in \mathcal{B}_{j(i-1)}^0} \varepsilon_k + c_{j(i)}^0 - c_{j(i-1)}^0 \right| > 0$$

hold with probability tending to 1 and $n \to \infty$.

We first consider (A.4). Notice that, for each $2 \leq i \neq j \leq n$, $\mathbb{E} d_i^\varepsilon = 0$, $\text{Var}\, d_i^\varepsilon = 2\sigma_n^2$ and

$$\text{Cov}(d_i^\varepsilon, d_j^\varepsilon) = \begin{cases} -\sigma_n^2, & \text{if } |i - j| = 1, \\ 0, & \text{otherwise.} \end{cases}$$

For $2 \leq i \leq n$, let $d_i^* \sim N(0, 2\sigma_n^2)$ be independent, so that

$$\begin{cases} \mathbb{E}(d_i^\varepsilon d_j^\varepsilon) \leq \mathbb{E}(d_i^* d_j^*), & \text{for all } 2 \leq i \neq j \leq n, \\ \mathbb{E}(d_i^\varepsilon)^2 = \mathbb{E}(d_i^*)^2, & \text{for all } 2 \leq i \leq n. \end{cases}$$

Then, by Slepian's inequality [see, e.g., Ledoux and Talagrand (1991)]

$$\mathbb{P}\Big\{\max_{i \in \mathcal{J}_0^c} |d_i^\varepsilon| \geq 4\lambda_{2,n}\Big\} \leq \mathbb{P}\Big\{\max_{i \in \mathcal{J}_0^c} |d_i^*| \geq 4\lambda_{2,n}\Big\}.$$

By Chernoff's bound for standard Gaussian variables, followed by the union bound

$$\mathbb{P}\Big\{\max_{i \in \mathcal{J}_0^c} |d_i^*| \geq 4\lambda_{2,n}\Big\} \leq 2\exp\Big\{-8\frac{\lambda_{2,n}^2}{\sigma_n^2} + \log |\mathcal{J}_0^c|\Big\},$$

which vanishes if condition 1 is satisfied.

In order to verify (A.5), it is sufficient to show that, with probability tending to 1 as $n \to \infty$,

$$\max_{i \in \mathcal{J}_0} \bigg| \frac{1}{b_{j(i)}^0} \sum_{k \in \mathcal{B}_{j(i)}^0} \varepsilon_k - \frac{1}{b_{j(i-1)}^0} \sum_{k \in \mathcal{B}_{j(i-1)}^0} \varepsilon_k + c_{j(i)}^0 - c_{j(i-1)}^0 \bigg| \leq \alpha_n,$$

where $\alpha_n = \min_{i \in \mathcal{J}_0} |d_i^0|$. By the triangle inequality, it is enough to show that

(A.6) $$\max_{i \in \mathcal{J}_0} \bigg| \frac{1}{b_{j(i)}^0} \sum_{k \in \mathcal{B}_{j(i)}^0} \varepsilon_k - \frac{1}{b_{j(i-1)}^0} \sum_{k \in \mathcal{B}_{j(i-1)}^0} \varepsilon_k \bigg| \leq \alpha_n/2$$

and

(A.7) $$\max_{i \in \mathcal{J}_0} |c_{j(i)}^0 - c_{j(i-1)}^0| \leq \alpha_n/2.$$

The previous inequality is implied by the last inequality in condition 2 in virtue of the bound

$$\max_{i \in \mathcal{J}_0} |c_{j(i)}^0 - c_{j(i-1)}^0| \leq 2\lambda_{2,n} \frac{1}{b_{\min}^0}.$$

Next, we turn to (A.6). Set $X_i = \frac{1}{b_{j(i)}^0} \sum_{k \in \mathcal{B}_{j(i)}^0} \varepsilon_k - \frac{1}{b_{j(i-1)}^0} \sum_{k \in \mathcal{B}_{j(i-1)}^0} \varepsilon_k$, with $i \in \mathcal{J}^0$. Then, $\mathbb{E} X_i = 0$ for all $i$ and

$$\max_{i \in \mathcal{J}^0} \text{Var } X_i \leq 2\frac{\sigma_n^2}{b_{\min}^0}.$$

Therefore, letting $X_i^* \sim N(0, 2\frac{\sigma_n^2}{b_{\min}^0})$, $i \in \mathcal{J}^0$, be independent, we obtain, using standard Gaussian tail bounds,

$$\mathbb{P}\Big\{\max_{i \in \mathcal{J}^0} |X_i| \geq \frac{\alpha_n}{2}\Big\} \leq \mathbb{P}\Big\{\max_{i \in \mathcal{J}^0} |X_i^*| \geq \frac{\alpha_n}{2}\Big\} \leq 2\exp\Big\{-\frac{b_{\min}^0 \alpha_n^2}{16\sigma_n^2} + \log |\mathcal{J}_0|\Big\}.$$

Under condition 2, the above probability vanishes. This, combined with (A.7) shows that (A.5) holds with probability tending to 1 if condition 2 is verified. □

PROOF OF THEOREM 2.5. It is enough to show that the event

$$\mathcal{R}_{\lambda_{1,n}} \cap \mathcal{R}_{\lambda_{2,n}}$$

occurs in probability for $n \to \infty$. Because the conditions of Theorem 2.3 are assumed, $\lim_n \mathbb{P}\{\mathcal{R}_{\lambda_{2,n}}\} = 1$, which implies that we can restrict our analysis to the set $\mathcal{R}_{\lambda_{2,n}}$, where $\widehat{J} = J_0$ and $\widehat{\mathcal{B}}_j = \mathcal{B}_j^0$, for $1 \leq j \leq J_0$. Next, from Corollary 2.2, it is immediately verified that the fused-lasso solution is

$$\widehat{\mu}^{\mathrm{FL}} = \sum_{j=1}^{\widehat{J}} 1_{\widehat{\mathcal{B}}_j} \widehat{\nu}_j^T,$$

where $\widehat{\nu}_j^T = \mathrm{sgn}(\widehat{\nu}_j)(\widehat{\nu}_j = \lambda_{1,n})_+$ is the soft-thresholded version of $\widehat{\nu}_j$. Therefore, in order to verify the claim, one needs to consider the simpler lasso problem applied to the vector $\widehat{\nu}$. Inspecting the sub-gradient condition for this problem, and by arguments similar to the ones used above, it follows that $\lim_n \mathbb{P}(\mathcal{R}_{\lambda_{1,n}}) = 1$ obtains provided both

$$(\mathrm{A.8}) \qquad \max_{j \in \mathcal{K}_0^c} \left| \frac{1}{b_j^0} \sum_{i \in \mathcal{B}_j^0} \varepsilon_i + c_j \right| < \lambda_{1,n}$$

and

$$(\mathrm{A.9}) \qquad \max_{j \in \mathcal{K}_0} \left| \frac{1}{b_j^0} \sum_{i \in \mathcal{B}_j^0} \varepsilon_i + c_j - \lambda_{1,n} \right| < \rho_n$$

hold with probability tending to 1 as $n \to \infty$, where the quantities $c_j$ are given in Lemma 2.1. Letting $X_j = \frac{1}{b_j^0} \sum_{i \in \mathcal{B}_j^0} \varepsilon_i$, notice that $X_j \sim N(0, \frac{\sigma_n^2}{b_j^0})$ and that $(X_1, \ldots, X_{J_0})$ are independent. Then, a combination of the Chernoff's and the union bounds yields

$$\mathbb{P}\left\{ \max_{j \in \mathcal{JS}_0^c} \left| \frac{1}{b_j^0} \sum_{i \in \mathcal{B}_j^0} \varepsilon_i \right| \geq \frac{\lambda_{1,n}}{2} \right\} \leq \sum_{j \in \mathcal{JS}_0^c} \exp\left\{ -\frac{\lambda_{1,n}^2 b_j^0}{8\sigma_n^2} \right\}$$

$$\leq \exp\left\{ -\frac{\lambda_{1,n}^2 b_{\min}^0}{8\sigma_n^2} + \log |\mathcal{JS}_0^c| \right\}$$

and

$$\mathbb{P}\left\{ \max_{j \in \mathcal{JS}_0} \left| \frac{1}{b_j^0} \sum_{i \in \mathcal{B}_j^0} \varepsilon_i \right| \geq \frac{\rho_n}{3} \right\} \leq \sum_{j \in \mathcal{JS}_0} \exp\left\{ -\frac{\rho_n^2 b_j^0}{18\sigma_n^2} \right\} \leq \exp\left\{ -\frac{\rho_n^2 b_{\min}^0}{18\sigma_n^2} + \log |\mathcal{JS}_0| \right\},$$

which give large deviations bounds for the error sums in (A.8) and (A.9). Conditions 1 and 3 guarantee that the above probabilities vanish for $n \to \infty$. Thus, with the additional conditions 2 and 4, the inequalities (A.8) and (A.9) are verified in probability. □

PROOF OF PROPOSITION 2.6. The proof is virtually identical to the proof of Theorem 2.5, the main differences stemming from the facts that the bias terms $c_j = 0$ for all $1 \leq j \leq J_0$ and

$$\frac{1}{\sqrt{b_j^0}} \sum_{i \in \mathcal{B}_j^0} \varepsilon_i \sim N(0, \sigma_n^2).$$

We omit the details. □

PROOF OF THEOREM 2.7. Let $\widehat{\mu}^F$ be the fusion estimate using the penalty $\lambda_{2,n}$. Then, because of assumption (2.14), and with the specific choice of $\lambda_{2,n}$ and $\sigma_n^2$ given in the statement, it can be verified that the conditions of Theorem 2.3 are met. Thus, the event

$$\mathcal{F} = \{\widehat{J} = J^0\} \cap \{\widehat{\mathcal{B}}_j = \mathcal{B}_j^0, 1 \leq j \leq J^0\}$$

has probability arbitrarily close to 1, for all $n$ large enough. On this event $\mathcal{F}$, we next investigate the adaptive fused-lasso $\widehat{\mu}$. Because $\widehat{\mu}$ is the minimizer of (2.10), for any $\mu \in \mathbb{R}^n$,

$$\|\widehat{\mu}^{\mathrm{AF}} - \widehat{\mu}\|_2^2 + 2\sum_i \lambda_i |\widehat{\mu}_i| \leq \|\widehat{\mu}^{\mathrm{AF}} - \mu\|_2^2 + 2\sum_i \lambda_i |\mu_i|,$$

where $\widehat{\mu}^{\mathrm{AF}}$ and $\lambda$ are given in (2.9) and (2.11), respectively. Adding and subtracting $\mu_0$ inside both terms $\|\widehat{\mu}^{\mathrm{AF}} - \widehat{\mu}\|_2^2$ and $\|\widehat{\mu}^{\mathrm{AF}} - \mu\|_2^2$ yields

(A.10) $\quad \|\widehat{\mu} - \mu_0\|_2^2 \leq \|\mu - \mu_0\|_2^2 + 2\sum_i \lambda_i(|\mu_i| - |\widehat{\mu}_i|) + 2\langle \varepsilon^*, \widehat{\mu} - \mu \rangle,$

where, on $\mathcal{F}$, $\varepsilon^* = \widehat{\mu}^{\mathrm{AF}} - \mu^0 = \sum_{j=1}^{J_0} X_j 1_{\mathcal{B}_j^0}$, with $X_j \sim N(0, \frac{\sigma_n^2}{b_j^0})$ and $(X_1, \ldots, X_{J^0})$ independent. Next, consider the sub-event $\mathcal{A} \subseteq \mathcal{F}$ given by

$$\mathcal{A} = \{|\varepsilon_i^*| \leq \lambda_i, \text{ for each } i = 1, \ldots, n\}$$
$$= \{|X_j| \leq \lambda_{1,n}/\sqrt{b_j^0}, \text{ for each } j = 1, \ldots, J_0\}.$$

Then,

$$\mathbb{P}(\mathcal{A}) = \mathbb{P}\Big\{\max_j |\zeta_j| \leq \lambda_{1,n}\Big\},$$

where $(\zeta, \ldots, \zeta_{J_0})$ are i.i.d. $N(0, \sigma_n^2)$. Notice that because of the choice of $\lambda_{1,n}$, $\lim_n \mathbb{P}\mathcal{A} = 1$ by standard large deviation bounds for Gaussians (see also the proof of Theorem 2.3). Next, on $\mathcal{A}$, we have

$$(\text{A.11}) \qquad 2\langle \varepsilon^*, \widehat{\mu} - \mu \rangle \leq 2 \sum_{i \in \mathcal{S}(\mu)} \lambda_i |\widehat{\mu}_i - \mu_i| + 2 \sum_{i \notin \mathcal{S}(\mu)} \lambda_i |\widehat{\mu}_i|.$$

The decomposition

$$2 \sum_i \lambda_i (|\mu_i| - |\widehat{\mu}_i|) = 2 \sum_{i \in \mathcal{S}(\mu)} \lambda_i |\mu_i| - 2 \sum_{i \in \mathcal{S}(\mu)} \lambda_i |\widehat{\mu}_i| - 2 \sum_{i \notin \mathcal{S}(\mu)} \lambda_i |\widehat{\mu}_i|,$$

along with (A.11) and the triangle inequality, yields, on $\mathcal{A}$,

$$2 \sum_i \lambda_i (|\mu_i| - |\widehat{\mu}_i|) + 2\langle \varepsilon^*, \widehat{\mu} - \mu \rangle \leq 4 \sum_{i \in \mathcal{S}(\mu)} \lambda_i |\widehat{\mu}_i - \mu_i|.$$

The previous display and (A.10) lead to the inequality

$$(\text{A.12}) \qquad \|\widehat{\mu} - \mu_0\|_2^2 \leq \|\mu - \mu_0\|_2^2 + 4 \sum_{i \in \mathcal{S}(\mu)} \lambda_i |\widehat{\mu}_i - \mu_i|$$

valid on $\mathcal{A}$. Next, it is easy to see that

$$\sum_{i \in \mathcal{S}(\mu)} \lambda_i^2 = \sum_{j \in \mathcal{JS}(\mu)} b_j \lambda_i^2 \leq \lambda_{1,n}^2 \sum_{l \in \mathcal{JS}(\mu;\mu^0)} 1 = \lambda_{1,n}^2 |\mathcal{JS}(\mu;\mu^0)|,$$

and, in particular,

$$\sum_{i \in \mathcal{S}(\mu)} \lambda_i^2 = \lambda_1^2 |\mathcal{JS}(\mu^0)|,$$

if and only if $\mathcal{JS}(\mu) = \mathcal{JS}(\mu^0)$.

Therefore, by the Cauchy–Schwarz inequality, the second term on the right-hand side of (A.12) can be bounded on $\mathcal{A}$ as follows:

$$4 \sum_{i \in \mathcal{S}(\mu)} \lambda_i |\widehat{\mu}_i - \mu_i| \leq 4 \lambda_{1,n} \sqrt{|\mathcal{JS}(\mu;\mu^0)|} \|\widehat{\mu} - \mu\|_2.$$

Then, using the triangle inequality, (A.12) becomes

$$\|\widehat{\mu} - \mu_0\|_2^2 \leq \|\mu - \mu_0\|_2^2 + 4 \lambda_{1,n} \sqrt{|\mathcal{JS}(\mu;\mu^0)|} (\|\widehat{\mu} - \mu^0\|_2 + \|\mu_0 - \mu\|_2).$$

On $\mathcal{A}$, the same arguments used in the second part of the proof of Lemma 3.7 in van de Geer (2007) establish the inequality in the claim. Since $\lim_n \mathbb{P}(\mathcal{A}) = 1$, the first result follows. $\square$

PROOF OF THEOREM 3.1. Let $N(\delta, \mathcal{F}_n, \|\cdot\|_n)$ denote the $\delta$-covering number of the set $\mathcal{F}_n \subset \mathbb{R}^n$ with respect to the norm $\|\cdot\|_n$ and notice that, for any $C > 0$,

$$N(\delta, C\mathcal{F}_n, \|\cdot\|_n) = N\left(\frac{\delta}{C}, \mathcal{F}_n, \|\cdot\|_n\right).$$

Furthermore, observe that $\mathcal{C}_{\mathrm{TV}}(T_n) = T_n\mathcal{C}(1)$. By a theorem of Birman and Solomjak (1967) [see, e.g., Lorentz, Golitschek and Makovoz (1996), Theorem 6.1], the $\delta$-metric entropy of $\mathcal{C}_{\mathrm{TV}}(T_n)$ with respect to the $L^2(\mathbb{P}_n)$ norm is

$$C\frac{T_n}{\delta},$$

for some constant $C$ independent of $n$. Letting $\Psi(\delta) = \int_0^\delta \sqrt{C\frac{T_n}{\delta}} = \sqrt{T_n C \delta}$, the solution to

$$\sqrt{n}\delta_n^2 \gtrsim \Psi(\delta_n)$$

gives

$$\delta_n \gtrsim \frac{T_n^{1/3}}{n^{1/3}},$$

where the symbol $\gtrsim$ indicates inequality up to a universal constant. The result now follows from Theorem 3.4.1 of van der Vaart and Wellner (1996) (see also the discussion on pages 331 and 332 of the same reference). In order to establish (3.2), we use Lemma 4.3 in Loubes and van de Geer (2002) to get that the metric entropy of $\mathcal{C}_{\ell_1}(L_n)$ is

$$H(\delta, \mathcal{C}_{\ell_1}(L_n), \|\cdot\|_n) \leq C\frac{L_n^2}{n\delta^2}\left(\log n + \log \frac{L_n}{\sqrt{n}\delta}\right)$$

for some constant $C$ independent of $n$. Notice that the entropy integral of $\sqrt{H(\delta, \mathcal{C}_{\ell_1}(L_n), \|\cdot\|_n)}$ diverges on any neighborhood of 0. By Theorem 9.1 in van de Geer (2000), the rate of consistency $\delta_n$ for $\widehat{\mu}^L$ with respect to the norm $\|\cdot\|_n$ is given by the solution to

(A.13) $$\sqrt{n}\delta_n^2 \gtrsim \Psi(\delta_n),$$

where

$$\Psi(\delta_n) \geq \int_{A\delta_n^2}^{\delta_n} \sqrt{H(x, \mathcal{C}_{\ell_1}(L_n))}\,dx$$

with $A$ a constant independent of $n$. Equation (A.13) is satisfied for a sequence $\delta_n$ satisfying

$$\sqrt{n}\delta_n^2 \gtrsim \frac{L_n\sqrt{\log n}}{\sqrt{n}}\log 1/\delta_n,$$

which gives the rate (3.2). $\square$

PROOF OF THEOREM 3.2. Let $H(\delta_n, \mathcal{C}(S_n, J_n), \|\cdot\|_n)$ denote the metric entropy of $\mathcal{C}(S_n, J_n)$ with respect to the norm $\|\cdot\|_n$. By Lemma A.2 and assumption (C2), for $\delta_n < 1$, the equation

$$\sqrt{n}\delta_n^2 \gtrsim \int_0^{\delta_n} \sqrt{\log H(x, \mathcal{C}(S_n, J_n), \|\cdot\|_n)}\, dx$$

leads to

$$\delta_n \gtrsim \sqrt{\frac{S_n}{n}} \log \sqrt{\frac{1}{\delta_n}} + o(1),$$

because $S_n = o(\frac{n}{\log n})$ and $j_n \leq s_n$. The sequence $\delta_n = \sqrt{\frac{J_n}{n}}$ satisfies the conditions of Theorem 3.4.1 of van der Vaart and Wellner (1996), thus proving (3.5). The second claim in the theorem is proved similarly, where the left-hand side of (A.14) in Lemma A.2 is now bounded by $C_{1,n}$ only. □

LEMMA A.2. *For the distance induced by the norm* $\|x\|_n = \frac{1}{\sqrt{n}}\sqrt{\sum_{i=1}^n x_i^2}$, *the metric entropy of* $\mathcal{C}(S_n, J_n)$ *satisfies*

(A.14) $$H(\delta, \mathcal{C}(S_n, J_n), \|\cdot\|_n) \leq C_{1,n} + C_{2,n},$$

*where*

$$C_{1,n} = \frac{J_n}{S_n}\log R_n + J_n\left(\log\frac{\sqrt{n}}{\delta} + \frac{1}{2}\log S_n\right) + S_n \log(S_n + J_n - 1)$$

*and*

$$C_{2,n} = \log S_n + S_n \log n.$$

PROOF. For fixed $\delta > 0$, we will construct an $\delta$-grid of $\mathcal{C}(S_n, J_n)$ based on the Euclidean distance. For every choice of $S_n$ nonzero entries of $\mu$, we regard $\mu$ as a vector in $\mathbb{R}^{S_n}$ which is block-wise constant with $J_n$ blocks. Then, there exist $J_n$ positive integer numbers $d_1, \ldots, d_{J_n}$ such that $\sum_l d_l = S_n$ and one can think of $\mu$ as the concatenation of $J_n$ vectors $\mu_1, \ldots, \mu_{J_n}$ each having constant entries, where $\mu_l \in \mathbb{R}^{d_l}$, $l = 1, \ldots, J_n$. Each $\mu_l$ can be any point along the main diagonal of the $d_l$-dimensional cube center at 0 with edge length $R_n^{1/S_n}$ and volume $R_n^{d_l/S_n}$. The length of the main diagonal of each such cube is $R_n^{1/S_n}\sqrt{d_l}$. Therefore, for any specific choice of $S_n$ nonzero coordinates, the slice in the corresponding $S_n$-dimensional cube centered at 0 and with edge length $R_n^{1/S_n}$ consisting of the set of vectors in $\mathcal{B}_n$ with discontinuity profile $(d_1, \ldots, d_{J_n})$ is the set

$$\mathcal{R}_n = \prod_{l=1}^{J_n} \ell(R_n^{1/S_n}, d_l),$$

where $\ell(R, d_l)$ denotes the closed line segment in $\mathbb{R}^{S_n}$ between the points $\pi_{d_l}(\mathbf{1}R)$ and $\pi_{d_l}(-\mathbf{1}R)$, where $\mathbf{1}$ is the $S_n$-dimensional vector with coordinates all equal to 1 and $\pi_{d_l}$ the function from $\mathbb{R}^{S_n}$ onto $\mathbb{R}^{S_n}$ given by $\pi_{d_l}(x) = y$ with $y_i = 0$ for $i \leq \sum_{j=1}^{l} d_l - 1$ or $i \geq \sum_{j=1}^{l+1} d_l$ and $y_i = x_i$ otherwise. Notice that the length of each $\ell(R_n^{1/S_n}, d_l)$ is precisely $R_n^{1/S_n} \sqrt{d_l}$. If $J_n = S_n$, $\mathcal{R}_n$ is the $S_n$-dimensional cube centered at 0 with volume $R_n$, while if $J_n < k_n$ the set $\mathcal{R}_n$ is a hyper-rectangle (not full dimensional) which can be embedded as a hyper-rectangle in $\mathbb{R}^{J_n}$ centered at 0 and with edge lengths equal to the lengths of $\ell(R_n^{1/k_n}, d_l)$, for $l = 1, \ldots, J_n$. As a result, it is immediate to see that the volume of $\mathcal{R}_n$ can be calculated as

$$\prod_l R_n^{1/S_n} \sqrt{d_l} = R_n^{J_n/S_n} \prod_l \sqrt{d_l}.$$

Next, partition each of the $J_n$ perpendicular sides of $\mathcal{R}_n$ into intervals of length $\delta\sqrt{\frac{d_l}{S_n}}$, $l = 1, \ldots, J_n$. This gives a partition of $\mathcal{R}_n$ into smaller hyper-rectangle of edge lengths $\delta\sqrt{\frac{d_l}{S_n}}$, for $l = 1, \ldots, J_n$. Every point in $\mathcal{R}_n$ is within Euclidean distance $\delta$ from the center of one of the small hyper-rectangles, which therefore form an $\delta$-grid for $\mathcal{R}_n$. By a volume comparison, the cardinality of such a grid is

$$\frac{R_n^{J_n/S_n} \prod_l \sqrt{d_l}}{\prod_l \delta\sqrt{d_l/S_n}} = \left(R_n^{1/S_n} \frac{\sqrt{S_n}}{\delta}\right)^{J_n}.$$

For fixed $S_n$, the number of distinct block patterns with cardinality at most $J_n$ is equal to the the number of nonnegative solutions to $d_1 + d_2 + \cdots + d_{J_n} = S_n$, which can bounded as

$$\binom{S_n + J_n - 1}{J_n} \leq (S_n + J_n - 1)^{J_n}$$

[see, e.g., Stanley (2000)]. Thus, the logarithm of cardinality of this $\delta$-grid is

(A.15) $\quad \dfrac{J_n}{S_n} \log R_n + J_n \left(\log \dfrac{1}{\delta} + \dfrac{1}{2} \log S_n\right) + J_n \log(S_n + J_n - 1).$

Next, the number of subsets of $\{1, \ldots, n\}$ of size at most $S_n$ is

$$\sum_{i=1}^{S_n} \binom{n}{i} \leq S_n n^{S_n}.$$

Thus, the logarithm of the cardinality for an $\delta$ grid over $\mathcal{B}_n$ is bounded by (A.15) plus the quantity

$$\log S_n + J_n \log n.$$

The result for the $\|\cdot\|_n$ norms now follows by replacing $\delta$ with $\delta/\sqrt{n}$ in (A.15). $\square$

**Acknowledgments.** I am very thankful to Larry Wasserman for indispensable advice, and to the anonymous referees and Associate Editor for valuable comments and suggestions.

Department of Statistics
Carnegie Mellon University
Pittsburgh, Pennsylvania 15213-3890
USA
E-mail: arinaldo@stat.cmu.edu